\documentclass[aap,seceqn,MSNbibl,citesort,dvips]{arximspdf}
\usepackage[vtexbibl]{}

\doi{10.1214/10-AAP697}
\volume{21}
\issue{1}
\pubyear{2011}
\firstpage{312}
\lastpage{331}

\makeatletter

\newtheorem{theorem}{Theorem}[section]
\newtheorem{proposition}[theorem]{Proposition}
\newtheorem{lemma}[theorem]{Lemma}
\newtheorem{corollary}[theorem]{Corollary}

\def\eqref#1{(\ref{#1})}
\makeatother

\begin{document}
\begin{frontmatter}

\title{Random permutations with cycle weights}
\runtitle{Random permutations with cycle weights}

\begin{aug}
\author[A]{\fnms{Volker} \snm{Betz}\thanksref{t1}\ead[label=e1]{v.m.betz@warwick.ac.uk}\ead
[label=u1,url]{http://www.maths.warwick.ac.uk/\texttildelow betz/}},
\author[A]{\fnms{Daniel} \snm{Ueltschi}\thanksref{t2}\corref{}\ead
[label=e2]{daniel@ueltschi.org}\ead[label=u2,url]{http://www.ueltschi.org}}
\and
\author[C]{\fnms{Yvan} \snm{Velenik}\thanksref{t3}\ead[label=e3]{Yvan.Velenik@unige.ch}\ead
[label=u3,url]{http://www.unige.ch/math/folks/velenik}}
\thankstext{t1}{Supported by the EPSRC fellowship EP/D07181X/1.}
 \thankstext{t2}{Supported in
part by the Grant DMS-0601075 of the US National Science Foundation.}
 \thankstext{t3}{Supported in part by the Swiss National Science Foundation.}

\runauthor{V. Betz, D. Ueltschi and Y. Velenik}
\affiliation{University of Warwick, University of Warwick and
Universit\'{e} de Gen\'{e}ve}
\address[A]{V. Betz\\
D. Ueltschi\\
Department of Mathematics\\
University of Warwick\\
Coventry, CV4 7AL\\
UK\\
\printead{e1}\\
\phantom{E-mail: }\printead*{e2}\\
\printead{u1}\\
\phantom{URL: }\printead*{u2}} 
\address[C]{Y. Velenik\\
Section de Math\'ematiques\\
Universit\'e de Gen\`eve\\
1211 Gen\`eve 4\\
Switzerland\\
\printead{e3}\\
\printead{u3}}
\end{aug}

\received{\smonth{8} \syear{2009}}
\revised{\smonth{2} \syear{2010}}

%
\begin{abstract}
We study the distribution of cycle lengths in models of nonuniform
random permutations with cycle weights.
We identify several regimes. Depending on the weights, the length of
typical cycles grows like the total
number $n$ of elements, or a fraction of $n$ or a logarithmic power of $n$.
\end{abstract}

%
\setattribute{keyword}{AMS}{AMS 2000 subject classification.}

\end{frontmatter}

\section{Introduction}\label{sec1}

We study the cycle distributions in models of weighted random
permutations. The probability of a permutation $\pi$ of $n$ elements
is defined by
%
\begin{equation}
\label{def prob}
P(\pi) = \frac1{h_n n!} \prod_{j\geq1} \theta_j^{r_j(\pi)},
\end{equation}
where $(\theta_1, \theta_2, \dots) \equiv{\bolds{\theta}}$ are real
nonnegative numbers, $r_j(\pi)$ denotes the number of $j$-cycles in
$\pi$ [we always have $\sum_j j r_j(\pi) = n$] and $h_n$ is the
normalization. We are mainly interested in the distribution of cycle
lengths in the limit $n \to\infty$ and in how these lengths depend on
the set of parameters ${\bolds{\theta}}$.

The probability $P$ is really a probability on sequences ${\mathbf r}=
(r_1,r_2,\dots)$ that satisfy $\sum_j j r_j = n$. It is well known
that ${\mathbf r}$ is the sequence of ``occupation numbers'' of a partition~${\bolds
{\lambda}}$
of $n$. That is, if ${\bolds{\lambda}}$ denotes the partition
$\lambda_1 \geq
\lambda_2 \geq\cdots$ with $\sum_i \lambda_i = n$, then $r_j$ is
the number of $\lambda_i$ that satisfies $\lambda_i = j$. Thus we are
really dealing with random partitions. The number of permutations that
are compatible with occupation numbers~${\mathbf r}$ is equal to
\[
\frac{n!}{\prod_{j\geq1} j^{r_j} r_j!}.
\]
It follows that the marginal of \eqref{def prob} on partitions is
given by
%
\begin{equation}
\label{distri partitions}
P({\bolds{\lambda}}) = \frac1{h_n} \prod_{j\geq1} \frac
1{r_j!} \biggl(
\frac1j \theta_j \biggr)^{r_j}.
\end{equation}
The formulas look simpler and more elegant for permutations than for
partitions and this is why we consider the former.

Random permutations with the uniform distribution have a compelling
history \cite{Gon,AT,DP,BP}. They are a special case of the present
setting, with $\theta_j \equiv1$.
The uniform distribution of random partitions has been studied, for
example, in \cite{EL,ST,Fri,VY}. They do not fit the present setting
because there are no parameters~${\bolds{\theta}}$ that make the~right-hand
side of \eqref{distri partitions} constant. Another distribution for
random partitions is~the Plancherel measure, where the probability of
${\bolds{\lambda}}$ is proportional to $\frac1{n!} (\dim{\bolds
{\lambda}})^2$; the
``dimension'' $\dim{\bolds{\lambda}}$ of a partition is defined as
the number
of Young tableaux in Young diagrams and it does not seem to have an
easy expression in terms of~${\mathbf r}$. Here again, we do not know
of any
direct relation between weighted random permutations and the Plancherel measure.

The present model was introduced in \cite{BU2} but variants of it have
been studied previously. The case of constant $\theta_j \equiv\theta
$ is known as the Ewens distribution. It appears in the study of
population dynamics in mathematical biology \cite{Ewe}; detailed
results about the number of cycles were obtained by Hansen \cite{Han}
and by Feng and Hoppe \cite{FH}. The distribution of cycle lengths was
considered by Lugo \cite{Lugo}. Another variant of this model involves
parameters $\theta_j \in\{0,1\}$, with finitely many 1's \cite
{Tim,Ben} or with parity dependence \cite{Lugo}.

Weighted random permutations also appear in the study of large systems
of quantum bosonic particles \cite{BU1,BU3}, where the parameters
${\bolds{\theta}}$ depend on such quantities as the temperature, the density
and the particle interactions. The $\theta_j$'s are thus forced upon
us and they do not necessarily take a simple form. This motivates the
present study where we only fix the asymptotic behavior of $\theta_j$
as $j \to\infty$.

The relevant random variables in our analysis are the lengths $\ell_i
= \ell_i(\pi)$ of the cycle containing the index $i = 1,\dots,n$.
These random variables are always identically distributed and obviously
not independent. Another relevant random variable is the number of
indices belonging to cycles of length between $a$ and $b$, $N_{a,b}(\pi
) = \#\{ i = 1,\dots,n  \dvtx  a\leq\ell_i(\pi)\leq b\}$.
It follows
from the exchangeability of $\ell_1,\dots,\ell_n$ that
%
\begin{equation}
\frac1n E(N_{a,b}) = P(\ell_1\in[a,b]).
\end{equation}
The properties of the distribution of $\ell_1$ that we derive below
can then be translated into properties of the expectation of $N_{a,b}$.

From a statistical mechanics point of view it is natural to introduce
the sequence ${\bolds{\alpha}}= (\alpha_1,\alpha_2,\dots)$ of parameters
such that ${\mathrm e}^{-\alpha_j} = \theta_j$. The model has an important
symmetry which is also a source of confusion, namely, the probability
of the permutation $\pi$ is left invariant under the transformation
%
\begin{equation}
\label{symmetry}
\alpha_j \mapsto\alpha_j + cj, \qquad h_n \mapsto{\mathrm e}^{-cn} h_n
\end{equation}
for any constant $c \in{\mathbb R}$. In particular, the case $\alpha
_j =
cj$ is identical to $\alpha_j \equiv0$, the case of uniform random
permutations.

The general results which we prove in this article rely on various
technical assumptions. To keep this Introduction simple, we only
describe the results in the particular but interesting case $\alpha_j
\sim j^\gamma$.
\begin{itemize}[$\bullet$]
\item The case $\gamma<0$ is a special case of the model studied
in~\cite{BU2} which is close to the uniform distribution.
\item In the case $\gamma=0$, that is,\ when $\theta_j \to\theta$
(the Ewens case, asymptotically), we find that $P(\ell_1 > sn) \to
(1-s)^\theta$.
Thus, almost all indices belong to cycles whose length is a fraction of
$n$. Precise statements and proofs can be found in Section \ref{sec Ewens}.
\item The case $0<\gamma<1$ is surprising. At first glance we might
expect smaller cycles than in the uniform case $\alpha_j \equiv0$.
However, we find that almost all indices belong to a single giant
cycle! The symmetry \eqref{symmetry} is indeed playing tricks on us.
In addition, we prove that the probability of the occurrence of a
single cycle of length $n$ is strictly positive and strictly less than
1. This is explained in detail in Section \ref{sec small alpha}.
\item The case $\gamma=1$ corresponds to uniform permutations because
of the symmetry \eqref{symmetry}.
\item When $\gamma>1$, the cycles become shorter and $\ell_1$ behaves
asymptotically as $(\frac1{\gamma-1} \log n)^{1/\gamma}$; see
Section \ref{sec alpha diverges quickly}.
\end{itemize}

Weighted random permutations clearly show a rich behavior and only a
little part has been uncovered so far. The case of negative parameters
$\alpha_j \asymp- j^\gamma$ remains to be explored and the future
will hopefully bring more results regarding concentration properties.

In the case of uniform permutations, it is known that the random
variables $r_k$ converge to independent Poisson random variables with
parameter $1/k$ in the limit $n\to\infty$~\cite{Gon,AT}. An open
problem is to understand how this generalizes to weighted random permutations.

\section{Asymptotic Ewens distribution}
\label{sec Ewens}

In the case of the uniform distribution, it is an easy exercise to show
that $P(\ell_1 = a) = 1/n$ for any $a = 1,\dots,n$. It follows that
$P(\ell_1 > s n) \to1-s$ for any $0 \leq s \leq1$. This
result was
extended to the case of small weights in \cite{BU2}. We consider here
parameters that are close to Ewens weights. A result similar to (a)
below has been recently derived by Lugo \cite{Lugo}.

\begin{theorem}
\label{thm Ewens}
Let $\theta\in{\mathbb R}_+$. We suppose that $\sum_{j=1}^\infty
\frac1j
|\theta_j - \theta| < \infty$ if $\theta\geq1$ or that $\sum
_{j=1}^\infty|\theta_j - \theta| < \infty$ if $\theta<1$.
\begin{longlist}[(a)]
\item[(a)] The distribution of $\ell_1$ satisfies, for $0\leq
s\leq1$,
%
\begin{equation}
\lim_{n\to\infty} P(\ell_1 > sn) = (1-s)^\theta.
\end{equation}
\item[(b)] The joint distribution of $\ell_1$ and $\ell_2$
satisfies, for $0\leq s,t\leq1$,
\begin{eqnarray}
&&\lim_{n\to\infty} P(\ell_1>sn,\ell_2>tn)\nonumber
\\[-8pt]\\[-8pt]
&&\qquad = \frac{\theta
}{1+\theta}
(1 - s - t)_+^{\theta+ 1} + \frac{1+\theta(s\vee t)}{1+\theta} (1 -
s\vee t)^\theta,\nonumber
\end{eqnarray}
\end{longlist}
where $f_+$ denotes the positive part of a function $f$.
\end{theorem}

Let us recall a few properties that are satisfied by the normalization
factors $h_n$.
Summing over the length $j$ of the cycle that contains 1 we find the
useful relation
%
\begin{equation}
\label{a useful relation indeed}
P(\ell_1 \in[a,b])= \frac1{n!h_n} \sum_{j=a}^b \frac
{(n-1)!}{(n-j)!} \theta_j (n-j)! h_{n-j} = \frac1n \sum_{j=a}^b
\theta_j \frac{h_{n-j}}{h_n}.
\end{equation}
Choosing $[a,b] = [1,n]$, we get
%
\begin{equation}
\label{useful consequence}
h_n = \frac1n \sum_{j=1}^n \theta_j h_{n-j}, \qquad h_0 = 1.
\end{equation}
Next, let $G_h(s) = \sum_{n\geq0} h_n s^n$ be the generating function
of the sequence $(h_n)$.
One can view a permutation as a combinatorial structure made of cycles.
It follows from standard combinatorics results that $G_h(s) =\break \exp\sum
_{j\geq1} \frac1j \theta_j s^j$. We also refer to \cite{BU2}
for a
direct proof of this formula. The first step in the proof of Theorem
\ref{thm Ewens} is to control the normalization $h_n$. Here, $(\theta
)_n = \theta(\theta+1) \cdots(\theta+n-1)$ denotes the ascending factorial.

\begin{proposition} \label{prop:ewens}
Under the assumptions of Theorem $\ref{thm Ewens}$, we have
\[
h_n = C({\bolds{\theta}}) \frac{(\theta)_n}{n!} \bigl(1 + o(1)\bigr) \qquad
\mbox{with }  C({\bolds{\theta}}) = \exp\sum_{j \geq1} \frac{1}{j}
(\theta_j - \theta) .
\]
\end{proposition}

\begin{pf}
We have
%
\begin{equation}
G_h(s) = \exp\biggl\{ \theta\sum_j \frac1j s^j + \sum_j \frac1j
(\theta_j - \theta) s^j \biggr\} = (1-s)^{-\theta} {\mathrm e}^{u(s)}
\end{equation}
with
%
\begin{equation}
u(s) = \sum_{j\geq1} \frac1j (\theta_j - \theta) s^j.
\end{equation}
Notice that $u(1) = \lim_{s \nearrow1} u(s)$ exists. Let $c_j$ be the
Taylor coefficients of ${\mathrm e}^{u(s)}$, that is, ${\mathrm
e}^{u(s)} = \sum c_j
s^j$. Then, by Leibniz' rule,
%
\begin{equation}
h_n = \frac{1}{n!} \frac{{\mathrm{d}}^n}{{\mathrm{d}}s^n} G_h(s)
\bigg\vert_{s=0} =
\frac{(\theta)_n}{n!} \sum_{k\geq0} d_{n,k} c_k
\end{equation}
with
%
\begin{equation}
d_{n,k} =
\cases{ \dfrac{n (n-1) \cdots(n-k+1)}{(\theta+n-1) \cdots(\theta
+n-k)}, &\quad\mbox{if } $k \leq n$,
\cr
0, &\quad\mbox{otherwise.}
}
\end{equation}
It is not hard to check that
%
\begin{equation}
\label{bounds for d}
d_{n,k} \leq
\cases{1, &\quad\mbox{if }$ \theta\geq1$,
 \cr
  \theta^{-1} k + 1, &\quad\mbox{if }$ \theta> 0$.
}
\end{equation}
Let $U(s) = \sum\frac1j |\theta_j - \theta| s^j$ and $C_j$ be the
Taylor coefficients of ${\mathrm e}^{U(s)}$. It is clear that $|c_j|
\leq C_j$
for all $j$. When $\theta\geq1$, the first bound of \eqref{bounds
for d} and the dominated convergence theorem imply
%
\begin{equation}
\label{ca converge !}
\lim_{n\to\infty} \sum_{k\geq0} d_{n,k} c_k = \sum
_{k\geq0} c_k
= {\mathrm e}^{u(1)} = C({\bolds{\theta}}).
\end{equation}
When $\theta<1$, the second bound of \eqref{bounds for d} gives
$d_{n,k} |c_k| \leq(\theta^{-1} k + 1) C_k$. The sequence $(k C_k)$
is absolutely convergent:
%
\begin{equation}
\qquad \sum k C_k = \frac{\mathrm{d}}{{\mathrm{d}}s} {\mathrm e}^{U(s)}
\bigg|_{s=1} = {\mathrm e}^{U(1)} U'(1)
= {\mathrm e}^{\sum(1/j) |\theta_j - \theta|} \sum|\theta_j -
\theta| <
\infty.
\end{equation}
We again obtain \eqref{ca converge !} by the dominated convergence theorem.
\end{pf}

\begin{pf*}{Proof of Theorem $\ref{thm Ewens}$}
We show that, for any $0<s<t<1$, we have
%
\begin{equation}
\lim_{n\to\infty} P(\ell_1 \in[sn,tn]) = (1-s)^\theta-
(1-t)^\theta.
\end{equation}
Using Proposition \ref{prop:ewens}, we have
%
\begin{equation}
\qquad P(\ell_1 \in[sn,tn]) = \frac1n \sum_{j=sn}^{tn} \theta_j \frac
{h_{n-j}}{h_n} = \frac\theta n \sum_{j=sn}^{tn} \frac{(\theta
)_{n-j}}{(n-j)!} \frac{n!}{(\theta)_n} \bigl(1 + o(1)\bigr).
\end{equation}
Here and throughout this article, when $a$ and $b$ are not integers we
use the convention
%
\begin{equation}
\sum_{j=a}^b f(j) = \sum_{j \in[a,b] \cap{\mathbb N}} f(j) = \sum
_{j =
\lceil a \rceil}^{\lfloor b \rfloor} f(j).
\end{equation}
We now use the identity
%
\begin{equation}
\label{Gamma identity}
(\theta)_n = \frac{\Gamma(n+\theta)}{\Gamma(\theta)}
\end{equation}
and the asymptotic
%
\begin{equation}
\label{Gamma asymptotic}
\frac{\Gamma(n+\theta)}{n!} = n^{\theta-1} \bigl(1 + o(1)\bigr).
\end{equation}
We get
%
\begin{equation}
P(\ell_1 \in[sn,tn]) = \frac\theta n \sum_{j=sn}^{tn} \biggl(1 - \frac
jn\biggr)^{\theta-1} \bigl(1+o(1)\bigr).
\end{equation}
As $n\to\infty$, the right-hand side converges to the Riemann
integral $\theta\int_s^t (1-\xi)^{\theta-1}\,
{\mathrm{d}}\xi$ and we obtain the first claim of Theorem \ref{thm Ewens}.

Let us now turn to the second claim. Let $1 \leq a \leq b
\leq n$ and
$1 \leq c \leq d \leq n$.
We get an expression for the joint probability of $\ell_1$ and $\ell
_2$ in a similar fashion as for \eqref{a useful relation indeed}. When
both indices belong to different cycles (noted $1 \not\sim2$), we have
%
\begin{equation}
\qquad\quad P ( \ell_1 \in[a,b], \ell_2 \in[c,d], 1 \not\sim2 ) =
\frac1{n! h_n} \mathop{\mathop{\sum_{j\in[a,b] }}_{k\in
[c,d]}}_{j+k \leq n} \mathop{\sum_{|c_1|=j}}_{|c_2|=k} \theta_j
\theta_k \sum_{\pi'} \prod_{\ell\geq
1} \theta_\ell^{r_\ell(\pi')}.
\end{equation}
Here $c_1$ and $c_2$ denote the cycles that contain 1 and 2,
respectively, and $\pi'$ denotes a permutation of the $n-j-k$ indices
that do not belong to $c_1$ or $c_2$. The number of cycles of length
$j$ that contain 1 but not 2 is $\frac{(n-2)!}{(n-1-j)!}$; given
$c_1$, the number\vspace*{-2pt} of cycles of length $k$ that contain 2 is $\frac
{(n-j-1)!}{(n-j-k)!}$. Since the sum over $\pi'$ gives $(n-j-k)!
h_{n-j-k}$, we get
%
\begin{equation}
P ( \ell_1 \in[a,b], \ell_2 \in[c,d], 1 \not\sim2 ) =
\frac1{n(n-1)} \mathop{\mathop{\sum_{j\in[a,b] }}_{k\in
[c,d]}}_{j+k \leq n} \theta
_j \theta_k \frac{h_{n-j-k}}{h_n}.
\end{equation}
When both indices belong to the same cycle one can first sum over the
length $j$ of the common cycle, then over $j-2$ indices other than 1, 2
and then over $j-1$ locations for 2. This gives $\frac{(n-2)!}{(n-j)!}
(j-1)$ possibilities. The sum over permutations on remaining indices
gives $(n-j)! h_{n-j}$. The result is
\begin{eqnarray}
P ( \ell_1 \in[a,b], \ell_2 \in[c,d] ) &=& \frac
1{n(n-1)} \mathop{\mathop{\sum_{j\in[a,b] }}_{k\in[c,d]}}_{j+k
\leq n} \theta_j
\theta_k \frac{h_{n-j-k}}{h_n}\nonumber
\\[-8pt]\\[-8pt]
&&{}+ \frac1{n(n-1)} \sum_{j\in
[a,b]\cap[c,d]} (j-1) \theta_j \frac{h_{n-j}}{h_n}.\nonumber
\end{eqnarray}
Let $\varepsilon>0$ and set $a=sn$, $c=tn$ and $b=d=n$. We assume,
without loss of generality, that $1\geq s\geq t\geq0$.
Using the above
expression, Proposition \ref{prop:ewens} and equations \eqref{Gamma
identity} and \eqref{Gamma asymptotic}, we deduce that, for $n$ large,
%
\begin{eqnarray}
&&P ( \ell_1 \geq sn, \ell_2 \geq tn )\nonumber
\\
&&\qquad= \frac
{\theta
^2}{n^2} \mathop{\sum_{j\geq sn, k\geq tn}}_{j+k \leq
(1-\varepsilon)n}
\biggl( 1 - \frac{j+k}{n} \biggr)^{\theta-1}   \bigl(1+o_\varepsilon(1)\bigr)
\\
&&\qquad\quad{}+ \frac\theta{n^2} \sum_{sn \leq j \leq(1-\varepsilon)n}
(j-1)
\biggl(1-\frac{j}{n}\biggr)^{\theta-1}   \bigl(1+o_\varepsilon(1)\bigr)  + O(\varepsilon).\nonumber
\end{eqnarray}
Taking first the limit $n\to\infty$ and then the limit $\varepsilon\to
0$, the right-hand side of the latter expression is seen to converge to
%
\begin{equation}
1_{\{s+t\leq1\}}\theta^2 \int_{s+t}^1 (\xi-s-t) (1-\xi
)^{\theta
-1} \,{\mathrm{d}}\xi+ \theta\int_s^1 \xi(1-\xi)^{\theta-1}
\,{\mathrm{d}}\xi
\end{equation}
and the second claim of Theorem \ref{thm Ewens} follows.
\end{pf*}

\section{Slowly diverging parameters}
\label{sec small alpha}

This section is devoted to parameters $\alpha_j$ that grow slowly to
$+\infty$. The typical case is $\alpha_j = j^\gamma$ with $0<\gamma
<1$ but our conditions allow more general sequences. As mentioned in
the \hyperref[sec1]{Introduction}, the system displays a surprising behavior: almost all
indices belong to a single giant cycle.

\begin{theorem}
\label{thm small alpha}
We assume that $0 < \frac{\theta_{n-j} \theta_j}{\theta_n}
\leq
c_j$ for all $n$ and for $j = 1,\dots, \frac n2$, with constants $c_j$
that satisfy $\sum_{j\geq1} \frac{c_j}j < \infty$. Then
\[
\lim_{m\to\infty} \lim_{n \to\infty} P(\ell_1 > n-m) = 1.
\]
\end{theorem}

It may be worth recalling that in this article $n$ always denotes the
number of elements and that $P$ depends on $n$. The proof of this
theorem can be found later in this section.
In the case $\alpha_j = j^\gamma$ we have
%
\begin{equation}
\frac{\theta_{n-j} \theta_j}{\theta_n} = {\mathrm e}^{-n^\gamma[
(1 -  j/n)^\gamma+ ( j/n)^\gamma- 1]} \approx
\cases{ {\mathrm e}^{-j^\gamma}, &\quad\mbox{if }$ j \ll n$,
\cr
{\mathrm e}^{- c n^\gamma}, &\quad\mbox{if }$ j = sn$,
}
\end{equation}
where the constant in the last equation is $c = (1-s)^\gamma+ s^\gamma
- 1$. It is positive for $0 < \gamma< 1$ and the condition of the
theorem is fulfilled. Another interesting example is $\theta_j =
j^{-\gamma}$ with $\gamma>0$, where we can choose $c_j = 2j^{-\gamma}$.

Let us understand why parameters $\alpha_j = j^\gamma$ favor longer
and longer cycles when $\gamma< 1$. The heuristics are actually
provided by statistical mechanics, namely, we can write the probability
$P(\pi)$ as a Gibbs distribution $\frac1Z {\mathrm e}^{-H(\pi)}$ with
``Hamiltonian'' $H(\pi) = \sum_{i=1}^n \frac{\alpha_{\ell_i(\pi
)}}{\ell_i(\pi)}$. Thus, an ``energy'' $\frac{\alpha_j}j =
j^{\gamma-1}$ is associated  with each index $i$ that belongs to a cycle of
length $j$.  Indices in longer cycles have lower energy so they are
favored. This discussion also provides an illustration for the symmetry
\eqref{symmetry}; it amounts to shifting the Hamiltonian by a constant
and this does not affect the Gibbs distribution.

We can state a more precise result than Theorem \ref{thm small alpha}
if we make the additional assumption that $\frac{\theta_{n+1}}{\theta
_n}$ converges to 1 as $n\to\infty$. This condition is easy to check
when $\alpha_j = j^\gamma$, $0<\gamma<1$ or when $\alpha_j = \gamma
\log j$, $\gamma>0$.

\begin{theorem}
\label{thm length 1}
Suppose that the assumptions of Theorem $\ref{thm small alpha}$ hold
true. In addition, we suppose that $\frac{\theta_{n+1}}{\theta_n}$
converges to $1$ as $n\to\infty$. Then $\sum_j h_j < \infty$, and
for any fixed $m\geq0$,
\[
\lim_{n\to\infty} P(\ell_1= n-m) = \frac{h_m}{\sum_{j\geq0} h_j}.
\]
\end{theorem}

Theorem \ref{thm length 1} shows in particular that a single cycle of
length $n$ occurs with probability $1 / \sum_j h_j$, but that finite
cycles may be present as well.

This theorem is proved at the end of the section. We first obtain
estimates for~$h_n$.

\begin{proposition}
\label{prop small alpha}
Under the assumptions of Theorem $\ref{thm small alpha}$ there exists
a constant $B$ such that, for all~$n\geq1$,
\[
1 \leq\frac{n h_n}{\theta_n} \leq B.
\]
The constant $B$ depends on $\{c_j\}$ only.
\end{proposition}

\begin{pf}
The lower bound follows obviously from \eqref{useful consequence} but
the upper bound requires some work.
Let $a_n = \frac{n h_n}{\theta_n}$.
The relation \eqref{useful consequence} can be written as
%
\begin{equation}
a_n = 1 + \sum_{j=1}^{n-1} \frac1j \frac{\theta_{n-j} \theta
_j}{\theta_n} a_j.
\end{equation}
We can rewrite this relation as
%
\begin{equation}
\quad \ \quad a_n =
\cases{ 1 + \displaystyle\sum_{j=1}^{{(n-1)}/2} \dfrac{\theta_{n-j}
\theta_j}{\theta_n} \biggl( \frac{a_j}j + \frac{a_{n-j}}{n-j}
\biggr), &\  \mbox{if $n$ is odd,}
\cr
1 + \displaystyle\sum_{j=1}^{ n/2 - 1} \dfrac
{\theta_{n-j} \theta_j}{\theta_n} \biggl( \frac{a_j}j + \frac
{a_{n-j}}{n-j} \biggr) + \frac{2 \theta_{n/2}^2}{n \theta_n} a_{n/2},
&\   \mbox{if $n$ is even.}
}
\end{equation}
We define the sequence $(b_n)$ by the recursion equation
%
\begin{equation}
b_n = 1 + \sum_{j=1}^{n/2} c_j \biggl( \frac{b_j}j + \frac
{b_{n-j}}{n-j} \biggr).
\end{equation}
It is clear that $a_n \leq b_n$ for all $n$. Next, let $m$ be a number
such that
%
\begin{equation}
\label{def m}
\frac2n \sum_{j=1}^{n/2} c_j + \sum_{j > m/2} \frac{c_j}j \leq
\frac12
\end{equation}
for all $n \geq m$. Such an $m$ exists because $(c_j/j)$ is summable
and the first term of the above equation is less than $\frac2{\sqrt n}
\sum_{j=1}^{\sqrt n} \frac{c_j}j + \sum_{j>\sqrt n} \frac{c_j}j =
o(1)$. We set
%
\begin{equation}
B = 2 \max_{1 \leq j \leq m} b_j.
\end{equation}
Notice that $B$ depends on the $c_j$s but not on the $\theta_j$s.
Finally, we introduce another sequence $(b_n')$ defined by
%
\begin{equation}
b_n' =
\cases{ b_n, &\quad\mbox{if }$ n \leq m$, \cr
 1 + \displaystyle\sum_{j=1}^{n/2} c_j
\biggl( \dfrac{b_j'}j + \dfrac{2B}n \biggr), &\quad\mbox{if }$ n>m$.
}
\end{equation}
It is clear that $b_n' \leq\frac12 B$ for $n \leq m$; we
now show by
induction that $b_n' \leq B$ for all~$n$. We have
\begin{eqnarray}
b_n' - b_m' &=& \frac{2B}n \sum_{j=1}^{n/2} c_j - \sum_{j=1}^{m/2}
c_j \frac{b_{m-j}'}{m-j} + \sum_{j =  m/2 + 1}^{n/2} c_j \frac
{b_j'}j\nonumber
\\[-8pt]\\[-8pt]
&\leq& \Biggl( \frac2n \sum_{j=1}^{n/2} c_j + \sum_{j > m/2}
\frac
{c_j}j \Biggr) B.\nonumber
\end{eqnarray}
This is less than $\frac12 B$ by definition \eqref{def m} of $m$.
Since $b_m' \leq\frac12 B$, we find that $b_n' \leq B$ for
all $n$.
The final step is to see that $b_n \leq b_n'$. This is clear when $n
\leq m$ and we get it by induction when $n>m$:
\begin{eqnarray}
 b_{n+1} &=& 1 + \sum_{j=1}^{n/2} c_j \biggl( \frac{b_j}j + \frac
{b_{n-j+1}}{n-j+1} \biggr)\nonumber
\\[-8pt]\\[-8pt]
&\leq& 1 + \sum_{j=1}^{n/2} c_j \biggl(
\frac
{b_j'}j + \frac{2B}{n+1} \biggr) = b_{n+1}'.\nonumber
\end{eqnarray}
We have shown that $a_n \leq b_n \leq b_n' \leq B$ for
all $n$.
\end{pf}

\begin{pf*}{Proof of Theorem \protect\ref{thm small alpha}}
Using Proposition \ref{prop small alpha} we get
\begin{eqnarray}
P(\ell_1 \leq n-m) &=& \frac1n \sum_{j=m}^{n-1} \theta_{n-j}
\frac
{h_j}{h_n} \leq B \sum_{j=m}^{n-1} \frac1j \frac{\theta_{n-j}
\theta_j}{\theta_n}\nonumber
\\[-8pt]\\[-8pt]
&\leq& B \sum_{j=m}^{n/2} \frac{c_j}j + B
\sum
_{j=n/2}^{n-1} \frac{c_{n-j}}j.\nonumber
\end{eqnarray}
The last term goes to zero as $n\to\infty$. The first term goes to
zero as $n\to\infty$ and $m\to\infty$.
\end{pf*}

\begin{pf*}{Proof of Theorem \ref{thm length 1}}
From equation\ \eqref{a useful relation indeed}
%
\begin{equation}
\label{thanks to the referee}
P(\ell_1 = n-m) = \frac1n \theta_{n-m} \frac{h_m}{h_n} = \frac
{\theta_{n-m}}{\theta_n} \frac{\theta_n}{n h_n} h_m.
\end{equation}
Further, \eqref{useful consequence} can be written as
%
\begin{equation}
\label{ca se domine}
\frac{n h_n}{\theta_n} = \sum_{j=0}^{n/2} \biggl( \frac{\theta
_{n-j} h_j }{\theta_n} + \frac{\theta_j h_{n-j}}{\theta_n} \biggr).
\end{equation}
This is actually correct for odd $n$ only; there is an unimportant
correction for even $n$ coming from $j=n/2$. Since $h_j \leq B
\frac
{\theta_j}j$ (Proposition \ref{prop small alpha}), the summand is
less than $B c_j (\frac1j + \frac1{n-j}) \leq2B \frac{c_j}j$. For
each $j$, and as $n\to\infty$, we have $\frac{\theta_{n-j}}{\theta
_n} \to1$ and $\frac{\theta_j h_{n-j}}{\theta_n} \leq B \frac
{c_j}{n-j} \to0$. The right-hand side of \eqref{ca se domine} then
converges to $\sum_j h_j$ by dominated convergence. We can now take
the limit $n\to\infty$ in \eqref{thanks to the referee} and we
indeed obtain the claim.
\end{pf*}

\section{Quickly diverging parameters}
\label{sec alpha diverges quickly}

Here we treat parameters $\theta_j = {\mathrm e}^{-\alpha_j}$ with
$\alpha_j$
diverging quickly, or equivalently
$\theta_j$ decaying quickly. More precisely, we shall make the
following two assumptions: for some $M>0$, all $k\geq1$ and two
coprime numbers $j_1,j_2 \geq4$,
%
\begin{equation}\label{funny assumpt}
0 \leq\theta_k \leq\frac{{\mathrm e}^{M k}}{k!},\qquad
\theta_{j_1}>0, \qquad \theta_{j_2}>0.
\end{equation}
It is necessary to impose some kind of aperiodicity condition on the
set of indices corresponding to nonvanishing coefficients $\theta_j$.
This prevents us from prescribing, for example,  permutations with
only even lengths of cycles. In this case we have $h_n=0$ for all odd
$n$, as can be easily seen from the recursion~\eqref{useful
consequence}; Proposition~\ref{h_n large coeff} below would fail.

Our assumptions allow us to get the asymptotics of $h_n$ using the
saddle point method. We write down the steps explicitly in order to
keep the article self-contained. A slightly shorter path would be to
prove that our assumptions imply that ${\mathrm e}^{f}$, with $f(z) =
\sum
_{j=0}^\infty\theta_j z^j$, is ``Hayman admissible'' and to use
standard results \cite{FS}.  Hayman admissibility  is
implicitly derived in our proof.

We describe general results in Section~\ref{ssec-general},
relegating proofs to Section~\ref{ssec-proofs}. The general results
turn out to be somewhat abstract so we use them to study the
particularly interesting class $\alpha_j=j^\gamma$, $\gamma>1$, in
Section~\ref{ssec-example}.

\subsection{Main properties}\label{ssec-general}
We now describe three general theorems about cycle lengths. In all
theorems conditions \eqref{funny assumpt} are silently assumed.
The first statement concerns the absence of macroscopic cycles.

\begin{theorem} \label{no large cycles}
For arbitrarily small $\delta> 0$ and arbitrarily large $k>0$, there
exists $n_0=n_0(\delta,k)$ such that
\[
P \Bigl( \max_{1\leq i\leq n}\ell_i \geq\delta n
\Bigr) \leq n^{-k}
\]
for all $n\geq n_0$.
\end{theorem}

More precise information about typical cycle lengths can be extracted
from the following result. Let $r_n>0$ be defined by the equation
%
\begin{equation}
\label{def rn}
\sum_{j\geq1} \theta_j r_n^j = n.
\end{equation}
That such $r_n$ exists uniquely is immediate.

\begin{theorem} \label{abstract large coeff}
Let $a(n), b(n)$ be such that
\[
\lim_{n \to\infty} \frac1n\sum_{j=1}^{a(n)} \theta_j r_n^{j+1/2}
= 0, \qquad
\lim_{n \to\infty} \frac1n\sum_{j=b(n)}^{n} \theta_j r_n^{j+1/2}
= 0.
\]
Then
\[
\lim_{n \to\infty} P\bigl(\ell_1 \in[a(n),b(n)]\bigr) = 1.
\]
\end{theorem}

When the information about the coefficients $\theta_j$ is sufficiently
detailed, some control on $r_n$ is possible and Theorem~\ref{abstract
large coeff} can be used to obtain sharp results. This is exemplified
in Section~\ref{ssec-example} for the special case $\alpha
_j=\alpha(j)=j^\gamma$ with $\gamma>1$.
In such cases, the sum $\sum_{j=1}^{\infty} \theta_j r_n^{j+1/2}$
(whose value is $r_n^{1/2} n$) is dominated by the terms corresponding
to indices $j$ close to the solution $j_{\max}$ of the equation
$\alpha
'(j)=\log r_n$.

Finally, it is also possible to extract from Theorem~\ref{abstract
large coeff} a general result proving absence of small cycles.

\begin{theorem} \label{small cycles}
\[
\lim_{n \to\infty} P \biggl( \ell_1 \leq\frac{\log n}{\log
r_n} -
\frac34 \biggr) = 0.
\]
\end{theorem}

We shall see below that the proof of Theorem \ref{small cycles} is
straightforward; nonetheless, the result is quite strong. In the case where
only finitely many $\theta_j$ are nonzero, we find
$r_n \sim n^{1/{j_0}}$, where $j_0$ is the last index with nonzero
$\theta_j$. Thus $\log n / \log r_n \approx j_0$ and
we obtain the probability that $\ell_1 \leq j_0-1$ is zero.
It follows that almost all cycles have length $j_0$, a fact already
observed in \cite{Tim,Ben}.
On the other hand, if infinitely many $\theta_j$ are nonzero,
it is easy to see that $\log n / \log r_n$ diverges. Thus $\ell_1$
goes to infinity in probability. To summarize, the only way to force a
positive fraction of indices to lie in finite cycles is to forbid
infinite cycles altogether, in which case
typical cycles have the maximal length that is allowed.

\subsection{Proofs of the main properties}
\label{ssec-proofs}

We now prove Theorems \ref{no large cycles}--\ref{small cycles}. We
use the following elementary result, which is a consequence of the
first assumption in~\eqref{funny assumpt}.

\begin{lemma}\label{lem derivative bound}
Let $f(x) = \sum_{k=0}^\infty c_k x^k$ with Taylor coefficients that satisfy
$0 \leq c_k \leq{\mathrm e}^{M k} k^{-k}$ for some $M>0$ and
all $k \geq1$. Then
for all $\delta>0$ and all $x\geq0$, we have
\[
f'(x) \leq(1+\delta) {\mathrm e}^{M} f(x) + {\mathrm
e}^{M}/\delta.
\]
\end{lemma}

\begin{pf}
Let $k_0 = k_0(x) = \lfloor(1+\delta) {\mathrm e}^{M} x \rfloor$. We
decompose
\[
f'(x) = \sum_{k=1}^\infty c_k k x^{k-1} = \sum_{k=1}^{k_0} c_k k
x^{k-1} + R(x).
\]
By our assumptions,
\begin{eqnarray*}
R(x) & =& \sum_{k=k_0+1}^\infty c_k k x^{k-1} \leq{\mathrm e}^{M}
\sum_{k=k_0+1}^\infty\biggl( \frac{x {\mathrm e}^{M}}{k}
\biggr)^{k-1}
\\
& \leq&{\mathrm e}^{M} \sum_{k=k_0+1}^\infty\biggl( \frac
{1}{1+\delta}
\biggr)^k
\leq{\mathrm e}^{M} /\delta.
\end{eqnarray*}
On the other hand, for the terms up to $k_0$, we have $k \leq k_0
\leq
(1+\delta) x {\mathrm e}^{M}$ and thus
\[
\sum_{k=1}^{k_0} c_k k x^{k-1} \leq(1+\delta) {\mathrm e}^{M}
\sum
_{k=0}^{k_0} c_k x^{k} \leq(1+\delta) {\mathrm e}^{M} f(x).
\]
This completes the proof.
\end{pf}

Let us define the functions
\[
I_\beta(z) = \sum_{j=1}^\infty j^\beta\theta_j z^j
\]
for $\beta\in{\mathbb R}$.
$\phi(z) := I_{-1}(z)$ plays a special role since the generating
function of $(h_n)$ is given by
$G_h(z) = \exp(\phi(z))$. All $I_\beta$ are analytic by the first
assumption in \eqref{funny assumpt}, monotone increasing and positive
on $\{z>0\}$
together with all their derivatives and $I_{\beta+1}(z) = z I_\beta'(z)
$. Lemma~\ref{lem derivative bound} implies that for each $\beta>0$ there
exists $C$ such that for all $z\geq0$ we have
%
\begin{equation}\label{main assumpt}
I_\beta'(z) \leq C I_\beta(z).
\end{equation}
Recall that $r_n = I_0^{-1}(n)$, where $ I_0^{-1}$ denote the inverse function.

\begin{proposition} \label{h_n large coeff}
We have
\[
h_n = \frac{r_n^{-n}}{\sqrt{2 \pi I_1(r_n)}} {\mathrm e}^{\phi(r_n)}
\bigl(1 + o(1)\bigr).
\]
\end{proposition}

\begin{pf}
Condition \eqref{funny assumpt} on Taylor coefficients implies that
$I_0(z) < \tilde D\times\break \exp(C z)$. Then
%
\begin{equation}\label{r_n cond}
r_n \geq c \log n
\end{equation}
for some $c>0$. On the other hand, $r_n$ diverges more slowly than
$n^{1/4}$ since $I_0(x)$ diverges faster than $x^4$ by \eqref{funny assumpt}.

For the saddle point method, we use Cauchy's formula and we obtain
\begin{eqnarray}\label{eq:split}
h_n &=& \frac1{2\pi r^n} \int_{-\pi}^\pi{\mathrm e}^{\phi(r
{\mathrm e}^{{\mathrm i}\gamma }) - n{\mathrm i}\gamma}\, {\mathrm
{d}}\gamma\nonumber
\\[-8pt]\\[-8pt]
&=& \frac{{\mathrm e}^{\phi(r)}}{2\pi r^n} \biggl[ \int_{-\gamma
_0}^{\gamma
_0} {\mathrm e}^{\phi(r {\mathrm e}^{{\mathrm i}\gamma}) - \phi(r) -
n{\mathrm i} \gamma}\, {\mathrm{d}}\gamma+ 2 \int_{\gamma_0}^\pi
{\mathrm e}^{\phi(r {\mathrm e}^{{\mathrm i} \gamma}) - \phi(r) -
n{\mathrm i}\gamma}\, {\mathrm{d}}\gamma\biggr]\nonumber
\end{eqnarray}
for any $r>0$ and any $0<\gamma_0<\pi$. We choose the $r = r_n$
defined by equation\ \eqref{def rn} since it is the minimum point of
$r^{-n} {\mathrm e}^{\phi(r)}$ and $\gamma_0 = \gamma_0(n) =
r_n^{-(1+\delta
)}$ for some $0<\delta<1/2$.
The leading order of the first term above can be found by expanding
$\phi(z)-n\log z$ around $\gamma=0$. We have
%
\begin{equation}
\phi(r_n {\mathrm e}^{{\mathrm i}\gamma}) - \phi(r_n) - n {\mathrm
i}\gamma= \sum_{j\geq
1} \frac{\theta_j}j r_n^j ( {\mathrm e}^{{\mathrm i}j \gamma}
- 1 - {\mathrm i}j
\gamma).
\end{equation}
Expanding ${\mathrm e}^{{\mathrm i}j \gamma} - 1 - {\mathrm i}j \gamma
= -\frac12 j^2
\gamma^2 + R(j\gamma)$
with $|R(j\gamma)| \leq\frac1{3!} (j\gamma)^3$ we get
\begin{eqnarray}
\phi(r_n {\mathrm e}^{{\mathrm i}\gamma}) - \phi(r_n) - n {\mathrm
i}\gamma&=& -\frac12
\gamma^2 \sum_{j\geq1} j \theta_j r_n^j + A(\gamma)\nonumber
\\[-8pt]\\[-8pt]
& =& -\frac
{1}{2} \gamma^2 I_1(r_n) + A(\gamma)\nonumber
\end{eqnarray}
with
%
\begin{equation}
|A(\gamma)| \leq\frac{\gamma_0^3}{3!} \sum_{j\geq1} j^2
\theta_j
r^j_n = \frac{\gamma^2}{r_n^{1+\delta} 3!} I_2(r_n)
\end{equation}
for all $\gamma\leq\gamma_0$.
Now, by \eqref{main assumpt}, we have $I_2(r_n) \leq C r_n I_1(r_n)$.
Thus, as $n \to\infty$, the term $A(\gamma)$
is negligible compared to $\gamma^2 I_1(r_n)$ in the first integral,
which is therefore given by
\begin{eqnarray}\label{first integral}
\quad \int_{-\gamma_0}^{\gamma_0} {\mathrm e}^{-{1/2} \gamma^2
I_1(r_n) (1 + o(1))} \,{\mathrm{d}}\gamma&=&
\frac{1}{\sqrt{I_1(r_n)}} \int_{-\gamma_0 \sqrt{I_1(r_n)}}^{\gamma
_0 \sqrt{I_1(r_n)}}
{\mathrm e}^{-{(1/2)} \xi^2 (1 + o(1))} \,{\mathrm{d}}\xi\nonumber
\\[-8pt]\\[-8pt]
&=& \sqrt
{\frac{2 \pi
}{I_1(r_n)}}\bigl(1 + o(1)\bigr).\nonumber
\end{eqnarray}
The last equality is justified by the fact that
$\gamma_0(n) I_1(r_n) \geq r_n^{-1-\delta} I_0(r_n) \geq
r_n^{-2} n$,
which diverges as $n \to\infty$.

We now turn to the second term in \eqref{eq:split}. We want to show
that it is negligible and we estimate it by
replacing the integral by $\pi$ times the maximum of the
integrand. In view of \eqref{first integral} it is enough to show that
%
\begin{equation}\label{need to show}
\lim_{n \to\infty} \frac{1}{2} \log I_1(r_n) - \operatorname{Re}\bigl(\phi(r_n)
- \phi(r_n {\mathrm e}^{{\mathrm i}\gamma})\bigr) = -\infty
\end{equation}
for all $\gamma\in[\gamma_0,\pi]$. For the first term we have
$\log I_1(r_n) \leq\log(C r_n I_0(r_n)) \leq\tilde C \log
n$. For
the second term
we have
\begin{eqnarray}
\qquad \operatorname{Re}\bigl(\phi(r_n) - \phi(r_n {\mathrm e}^{{\mathrm i}\gamma})\bigr)
&=&
\sum_{j \geq1}
\frac{1}{j} \theta_j r_n^j \bigl(1 - \cos(\gamma j)\bigr)\nonumber
\\[-8pt]\\[-8pt]
&\geq&\frac{\theta_{j_1}}{j_1} r_n^{j_1} \bigl(1-\cos(\gamma j_1)\bigr) +
\frac{\theta_{j_2}}{j_2} r_n^{j_2} \bigl(1-\cos(\gamma j_2)\bigr),\nonumber
\end{eqnarray}
where $j_1$ and $j_2$ are picked according to \eqref{funny assumpt}.
The right-hand side is zero at $\gamma=0$ and it is strictly positive
when $\gamma\in(0,\pi]$ ($j_1$ and $j_2$ are coprime); so its
minimum is attained at $\gamma_0$ when $n$ is sufficiently large
(recall that $\gamma_0 \to0$ when $n\to\infty$). Expanding the
cosine, we get
%
\begin{equation}
\operatorname{Re}\bigl(\phi(r_n) - \phi(r_n {\mathrm e}^{{\mathrm i}\gamma})\bigr)
\geq c' r_n^4 \gamma
_0^2 = c' r_n^{2 - 2 \delta} \geq c c' (\log
n)^{2-2\delta}.
\end{equation}
This dominates the first term of \eqref{need to show} since $\delta
<1/2$ and the proof is complete.
\end{pf}

\begin{pf*}{Proof of Theorem~$\ref{no large cycles}$}
Clearly,
%
\begin{equation}
P\Bigl(\max_i \ell_i > \delta n\Bigr) \leq n P(\ell_1 > \delta n).
\end{equation}
We have $I_1(r_n) \leq C^2 r_n^2 \phi(r_n)$ by \eqref{main assumpt}
and thus
Proposition \ref{h_n large coeff} gives $h_n \geq C' r_n^{-n-1}$
for $n$
large enough. Since all the $h_{n-j}$'s are clearly bounded by some
$D>0$, we have by \eqref{a useful relation indeed}
\begin{eqnarray}
n P(\ell_1 > \delta n) &\leq& D r_n^{n+1} \sum_{j = \delta n}^n
\biggl( \frac{{\mathrm e}^{M}}{j} \biggr)^j
\leq D r_n^{n+1} n \biggl(\frac{{\mathrm e}^{M}}{\delta n}
\biggr)^{\delta n}\nonumber
\\[-8pt]\\[-8pt]
&\leq& D n \biggl(\frac{{\mathrm e}^{M} r_n^{2/\delta}}{\delta
n}\biggr)^{\delta n}.\nonumber
\end{eqnarray}
The statement is trivial [and seen directly from \eqref{a useful relation
indeed}] if only finitely many $\theta_j$ are nonzero; thus we may assume
there are infinitely many nonzero $\theta_j$. Then
$I_0(z)$ grows faster at infinity than any power of $z$ and $r_n$
diverges more slowly than any power of $n$. The last bracket is less
than 1 for $n$ large enough so that the right-hand side vanishes in the
limit $n \to\infty$.
\end{pf*}

In order to make more precise statements about the length of typical
cycles we need a better control over the terms appearing in \eqref{a
useful relation indeed}. By the previous result it suffices to consider
the case where $j$ is not too close to $n$.

\begin{proposition} \label{h_n ratios}
For each $\delta> 0$ there exists $C_\delta$ such that,
for all $n \in{\mathbb N}$ and all $j < (1-\delta) n$, we have
\[
\frac{h_{n-j}}{h_n} \leq C_\delta r_n^{j+1/2}.
\]
\end{proposition}

\begin{pf}By Proposition \ref{h_n large coeff} we have
%
\begin{eqnarray}
\frac{h_{n-j}}{h_n} & \approx& r_n^j \biggl( \frac{r_n}{r_{n-j}}
\biggr)^{n-j}
\biggl( \frac{I_1(r_n)}{I_1(r_{n-j})} \biggr)^{1/2} {\mathrm
e}^{\phi (r_{n-j})-\phi(r_n)}\nonumber
\\
&=& r_n^j \exp\bigl( - \bigl( \phi(r_n) - \phi(r_{n-j})\nonumber
\\
&&{}\hspace*{41pt}- (n-j) \bigl(\ln
(r_n) - \ln(r_{n-j})\bigr) \bigr) \bigr)
\biggl( \frac{I_1(r_n)}{I_1(r_{n-j})} \biggr)^{1/2}
\\
&=& r_n^j \exp\biggl( - \biggl( \phi(r_n) - \phi(r_{n-j})\nonumber
\\
&&{}\hspace*{45pt}- \phi
'(r_{n-j}) r_{n-j} \ln\biggl(\frac{r_n}{r_{n-j}}\biggr) \biggr) \biggr)
\biggl( \frac{I_1(r_n)}{I_1(r_{n-j})} \biggr)^{1/2}\nonumber
\end{eqnarray}
when both $n$ and $n-j$ are large.
Put $r_{n-j} = x$ and $r_n = x + u$. Since $n \mapsto r_n$ is
increasing, we have $u > 0$. The exponent above then has the form
\begin{eqnarray}
\quad&&\phi(x+u) - \phi(x) - x \phi'(x) \ln\biggl(\frac{x+u}{x}\biggr)\nonumber
\\[-8pt]\\[-8pt]
&&\qquad = \bigl(\phi(x+u)
- \phi(x) - \phi'(x) u\bigr) + \phi'(x) \biggl(u - x \ln\biggl(\frac{x+u}{x}\biggr)\biggr).\nonumber
\end{eqnarray}
The first bracket in the right-hand side is greater than $\frac12 u^2
\phi''(x)$ since all derivatives
of $\phi$ are positive on ${\mathbb R}^+$. The second bracket is always
positive. Thus, for all $n \in{\mathbb N}$ and all
$j \leq(1 - \delta) n$, there exists $C_\delta' > 0$ such that
%
\begin{equation}
\frac{h_{n-j}}{h_n} \leq C_\delta' r_n^j {\mathrm e}^{-
(1/2) (r_n - r_{n-j})^2 \phi''(r_{n-j})}
\biggl( \frac{I_1(r_n)}{I_1(r_{n-j})} \biggr)^{1/2}.
\end{equation}
By \eqref{main assumpt}, $I_1(x) = x I_0'(x) \leq C x I_0(x)$. We
also have
$I_0(x) \leq I_1(x)$. Since $I_0(r_n) = n$, we get
%
\begin{equation}
\frac{I_1(r_n)}{I_1(r_{n-j})} \leq C r_n \frac{n}{n-j} \leq
\frac
{C}{\delta} r_n.
\end{equation}
This proves the claim.
\end{pf}

\begin{pf*}{Proof of Theorem \ref{abstract large coeff}}
The claims follows immediately from \eqref{a useful relation indeed}
and Proposition~\ref{h_n ratios}.
\end{pf*}

\begin{pf*}{Proof of Theorem \ref{small cycles}}
Let $m = \log n / \log r_n - \frac34$. We use equation\ \eqref{a
useful relation indeed}, bounding $\theta_j$ by a constant and using
Proposition \ref{h_n ratios} for the ratio of normalization factors.
Since $r_n$ diverges, we have
%
\begin{equation}
P(\ell_1 \leq m) \leq\frac{C}{n} \sum_{j=1}^m r_n^{j +
1/2} =
\frac{C}{n} r_n^{3/2} \frac{r_n^m-1}{r_n-1} \leq\frac{C'}{n}
r_n^{m+1/2},
\end{equation}
if $n$ is large enough. The right-hand side is equal to $C' r_n^{-1/4}$
and it vanishes in the limit $n\to\infty$.
\end{pf*}

\subsection{An explicit example}
\label{ssec-example}

In this subsection we treat explicitly the case $\alpha_j = \alpha(j)
= j^\gamma$ with $\gamma> 1$ as an example of application of the
previous general results. We first observe that the assumptions~\eqref
{funny assumpt} are trivially satisfied so that the general results in
this section apply.

The main result of this subsection is that typical cycles are of
size\break
$(\frac1{\gamma-1}\log n)^{1/\gamma}$ to leading order.

\begin{theorem}\label{thm-logsize}
Let $\alpha_j = j^\gamma$, with $\gamma>1$. Then
%
\begin{equation}
\frac{\ell_1}{((1/({\gamma-1}))\log_n)^{1/\gamma}} \to1
\end{equation}
in probability.
\end{theorem}

Let us define
%
\begin{equation}
\Delta(j) = \alpha(j) - \alpha(j_{\max}) - (j-j_{\max}) \log r_n.
\end{equation}
The proof of Theorem~\ref{thm-logsize} follows from two simple
technical estimates.

\begin{lemma}\label{lem-simple}
Let $j_{\max}\in{\mathbb R}$ be such that $\alpha'(j_{\max}) = \log r_n$.
\begin{enumerate}[(a)]
\item[(a)] Assume that $\gamma\geq2$.
Then for all $j\geq1$, there exists $c=c(\gamma)>0$ such that
%
\begin{equation}
\Delta(j) \geq c \alpha''(j_{\max})(j-j_{\max})^2.
\end{equation}
(When $j\geq j_{\max}$, one can choose $c=\frac12$.)
\item[(b)] Assume that $\gamma\in(1,2)$.
Then, for all $1\leq j\leq2j_{\max}$, there exists
$c=c(\gamma)>0$ such that
\begin{equation}
\Delta(j) \geq c \alpha''(j_{\max})(j-j_{\max})^2.
\end{equation}
(When $j\leq j_{\max}$, one can choose $c=\frac12$.)
Moreover, for all $j>2j_{\max}$, there exists $c=c(\gamma)>0$ such that
%
\begin{equation}
\Delta(j) \geq c j^\gamma.
\end{equation}
\end{enumerate}
\end{lemma}

\begin{pf}
We start with the case $\gamma\geq2$.
First of all, since $j_{\max}=(\alpha')^{-1}(\log r_n)$, we have for any
$j>j_{\max}$
%
\begin{eqnarray}
\Delta(j) &=& \alpha(j) - \alpha(j_{\max}) - (j-j_{\max}) \log r_n\nonumber
\\
&=&
\alpha(j) - \alpha(j_{\max}) - (j-j_{\max}) \alpha'(j_{\max})
\\
&=& \int_{j_{\max}}^j {\mathrm{d}}s\int_{j_{\max}}^s \alpha''(t)
\,{\mathrm{d}}t
\geq\frac12 \alpha''(j_{\max}) (j-j_{\max})^2,\nonumber
\end{eqnarray}
since $\alpha''$ is an increasing function.
Similarly, we have for any $\frac12j_{\max}\leq j<j_{\max}$
\begin{eqnarray}
\Delta(j) &=& \int_j^{j_{\max}}{\mathrm{d}}s \int_s^{j_{\max}}
\alpha''(t) \,{\mathrm{d}}t\nonumber
\\
&\geq&\frac12 \alpha''\biggl(\frac12j_{\max}\biggr) (j-j_{\max})^2
\\
&=& 2^{1-\gamma} \alpha''(j_{\max}) (j-j_{\max})^2.\nonumber
\end{eqnarray}
Finally, for $0\leq j<\frac12j_{\max}$ we use
\begin{eqnarray}
\Delta(j) &=& \int_j^{j_{\max}}\,{\mathrm{d}}s \int_s^{j_{\max}}
\alpha''(t) \,{\mathrm{d}}t
\geq\int_{j_{\max}/2}^{j_{\max}}\,{\mathrm{d}}s \int
_s^{j_{\max}} \alpha''(t) \,{\mathrm{d}}t\nonumber
\\[-8pt]\\[-8pt]
&\geq&\frac12 \alpha''\biggl(\frac12j_{\max}\biggr) \frac14 j_{\max}^2
\geq2^{-\gamma-1} \alpha''(j_{\max}) (j-j_{\max})^2.\nonumber
\end{eqnarray}

Let us now turn to the case $\gamma\in(1,2)$. The proof is completely
similar. When $j\leq j_{\max}$ we use (observe that $\alpha''$
is a
decreasing function now)
%
\begin{equation}
\Delta(j) = \int_j^{j_{\max}} \,{\mathrm{d}}s\int_s^{j_{\max}}
\alpha''(t) \,{\mathrm{d}}t
\geq\frac12 \alpha''(j_{\max}) (j-j_{\max})^2.
\end{equation}
When $j_{\max}<j\leq2j_{\max}$ we use
\begin{eqnarray}
\Delta(j) &=& \int_{j_{\max}}^j\,{\mathrm{d}}s \int_{j_{\max}}^s
\alpha''(t) \,{\mathrm{d}}t
\geq\frac12 \alpha''(2j_{\max}) (j-j_{\max})^2\nonumber
\\[-8pt]\\[-8pt]
&=& 2^{\gamma-3} \alpha''(j_{\max}) (j-j_{\max})^2.\nonumber
\end{eqnarray}
Finally, when $j>2j_{\max}$ we have
\begin{eqnarray}
\Delta(j) &=& \int_{j_{\max}}^j\,{\mathrm{d}}s \int_{j_{\max}}^s
\alpha''(t) \,{\mathrm{d}}t
\geq\frac12 \alpha''(j) (j-j_{\max})^2\nonumber
\\[-8pt]\\[-8pt]
&\geq&\frac18 \alpha''(j) j^2
= \frac18\gamma(\gamma-1) j^\gamma.\nonumber
\end{eqnarray}
\upqed\end{pf}

\begin{corollary}\label{Cor-asympt-gamma}
For any $\gamma>1$, we have, as $n\to\infty$,
%
\begin{eqnarray}\label{eq-jmax}
j_{\max}&=& \biggl(\frac1{\gamma-1}\log n\biggr)^{1/\gamma} \bigl(1+o(1)\bigr),
\\
\label{eq-rn}
\log r_n &=& \alpha'(j_{\max}) = \gamma\biggl(\frac1{\gamma-1}\log
n\biggr)^{(\gamma-1)/\gamma} \bigl(1+o(1)\bigr),
\\\label{eq-max}
{\mathrm e}^{-\alpha(j_{\max})}r_n^{j_{\max}} &=& n^{1+o(1)}.
\end{eqnarray}
\end{corollary}

\begin{pf}
We start with the case $\gamma\geq2$. Using the previous lemma, it
immediately follows that
\begin{eqnarray}\label{eq-I0-ub}
I_0(r_n) &=& \sum_{j\geq1} \mathrm{e}^{-\alpha(j)} r_n^j \leq
\mathrm{e}^{-\alpha(j_{\max}
)}r_n^{j_{\max}} \sum_{j\geq1} \mathrm{e}^{-c\alpha''(j_{\max
})(j-j_{\max})^2}\nonumber
\\[-8pt]\\[-8pt]
&\leq& C_1  \mathrm{e}^{-\alpha(j_{\max})}r_n^{j_{\max}}.\nonumber
\end{eqnarray}
Since for $j<j_{\max}$, $\Delta(j)\leq\frac12\alpha
''(j_{\max}
)(j-j_{\max})^2$, we also have
%
\begin{equation}\label{eq-I0-lb}
I_0(r_n) \geq \mathrm{e}^{-\alpha(\lfloor j_{\max}\rfloor)}r_n^{\lfloor
j_{\max}
\rfloor} \geq \mathrm{e}^{-(1/2) \alpha''(j_{\max})}  \mathrm{e}^{-\alpha
(j_{\max}
)}r_n^{j_{\max}}.
\end{equation}
Using the relation $I_0(r_n) = n$, \eqref{eq-I0-ub} and~\eqref
{eq-I0-lb} immediately imply the claimed asymptotics.

Let us now turn to the case $\gamma\in(1,2)$.
The lemma implies that
\begin{eqnarray}
I_0(r_n) &=& \mathrm{e}^{-\alpha(j_{\max})}r_n^{j_{\max}} \sum_{j\geq1}
\mathrm{e}^{-\Delta(j)}\nonumber
\\[-8pt]\\[-8pt]
&\leq& C_2  \mathrm{e}^{-\alpha(j_{\max})}r_n^{j_{\max}} \biggl\{\alpha
''(j_{\max}
)^{-1/2} + \sum_{j>2j_{\max}} \mathrm{e}^{-c j^\gamma} \biggr\}.\nonumber
\end{eqnarray}
Since $j_{\max}\nearrow\infty$ as $n\to\infty$, we see that $\sum
_{j>2j_{\max}} \mathrm{e}^{-c j^\gamma} \ll\alpha''(j_{\max})^{-1/2}$ and thus
that, for large $n$,
%
\begin{equation}
I_0(r_n) \leq C_3  \alpha''(j_{\max})^{-1/2}   \mathrm{e}^{-\alpha
(j_{\max}
)}r_n^{j_{\max}}.
\end{equation}
As above, we also have
\begin{eqnarray}\label{eq-I0-lb-geq2}
I_0(r_n) &\geq& \mathrm{e}^{-\alpha(\lceil j_{\max}\rceil)}r_n^{\lceil
j_{\max}
\rceil} \geq \mathrm{e}^{-(1/2) \alpha''(j_{\max})}  \mathrm{e}^{-\alpha
(j_{\max}
)}r_n^{j_{\max}}\nonumber
\\[-8pt]\\[-8pt]
&\geq& C_4  \mathrm{e}^{-\alpha(j_{\max})}r_n^{j_{\max}}.\nonumber
\end{eqnarray}
The claimed asymptotics follow as before.
\end{pf}

\begin{pf*}{Proof of Theorem~\ref{thm-logsize}}
Let $\varepsilon>0$. It is sufficient to check that Theorem~\ref
{abstract large coeff} applies with $a(n) = (1-\varepsilon)j_{\max}$ and
$b(n)=(1+\varepsilon)j_{\max}$. It follows from~Lemma~\ref{lem-simple} and
Corollary~\ref{Cor-asympt-gamma} that
%
\begin{equation}
\frac1n \sum_{j=b(n)}^\infty \mathrm{e}^{-\alpha(j)} r_n^{j+1/2} \leq
n^{o(1)} \sum_{j=b(n)}^\infty \mathrm{e}^{-c \alpha''(j_{\max}) (j-j_{\max})^2},
\end{equation}
which goes to $0$ as $n\to\infty$, since
%
\begin{equation}
\mathrm{e}^{-c\alpha''(j_{\max}) (b(n)-j_{\max})^2} = n^{-c\varepsilon^2\gamma
(1+o(1))}.
\end{equation}
Similarly,
\begin{eqnarray}
\frac1n \sum_{j=1}^{a(n)} \mathrm{e}^{-\alpha(j)} r_n^{j+1/2} &\leq& n^{o(1)}
\sum_{j=1}^{a(n)} \mathrm{e}^{-c \alpha''(j_{\max}) (j-j_{\max})^2}\nonumber
\\[-8pt]\\[-8pt]
&\leq&
n^{o(1)} \mathrm{e}^{-c \alpha''(j_{\max})j_{\max}^2 \varepsilon^2},\nonumber
\end{eqnarray}
which again goes to $0$ as $n\to\infty$.
\end{pf*}

\section*{Acknowledgments}
We are indebted to the referee for several useful comments and
especially for suggesting the claim of Theorem~\ref{thm length 1}.
D. Ueltschi is grateful to Nick Ercolani and several members of the University
of Arizona for many discussions about the Plancherel measure. D. Ueltschi also
acknowledges the hospitality of the University of Geneva, ETH Z\"urich,
the Center of Theoretical Studies of Prague and the University of
Arizona where parts of this project were carried forward.


\printaddresses

\end{document}